\documentclass{lmcs}
\pdfoutput=1

\usepackage{lastpage}
\lmcsdoi{22}{2}{1}
\lmcsheading{}{\pageref{LastPage}}{}{}%
{Jul.~01,~2025}{Apr.~07,~2026}{}

\usepackage{mathtools}
\usepackage{amsmath}
\usepackage{amssymb}
\usepackage{amsthm}
\usepackage{graphicx} 
\usepackage{hyperref}
\usepackage[utf8]{inputenc}

\newcommand{\PV}{\text{PV}}
\newcommand{\abs}[1]{\lvert #1\rvert}
\renewcommand{\P}{\mathcal{P}}
\newcommand{\F}{\mathcal{F}}
\newcommand{\C}{\mathcal{C}}
\newcommand{\floor}[1]{\lfloor #1\rfloor}
\newcommand{\PrimeFactor}{\textsc{PrimeFactor}}
\newcommand{\NN}{\mathbb{N}}

\begin{document}

\title{Prime Factorization in Models of $\text{PV}_1$}
\author[O.~Ježil]{Ondřej Ježil}
\address{Faculty of Mathematics and Physics, Charles University}
\email{ondrej.jezil@email.cz}
\thanks{This work has been supported by Charles University Research Center program No.~UNCE/24/SCI/022,
the project SVV-2025-260837 and by the GA UK project No. 246223.} 

\begin{abstract}
Assuming that no family of polynomial-size Boolean circuits can factorize a constant fraction of all products of two $n$-bit primes, we show that the bounded arithmetic theory $\PV_1$, even when augmented by the sharply bounded choice scheme $BB(\Sigma^b_0)$, cannot prove that every number has some prime divisor. By the completeness theorem, it follows that under this assumption there is a model $M$ of $\PV_1$ that contains a nonstandard number $m$ which has no prime factorization.
\end{abstract}

\maketitle

\section{Introduction}

Bounded arithmetic is a collective name for a family of first-order theories, which are weak fragments of Peano arithmetic, with strong connections to complexity theory. Their axiomatization usually consists of some basic universal theory describing the recursive properties of the function and relation symbols in the language of the theory, which usually extends the language of Peano arithmetic $L_{PA}=\{0,1,+,\cdot, \leq\}$, and an induction or some other scheme for a class of formulas whose expressivity corresponds to some computational complexity class. For more details about the motivation behind the study of these theories, we refer the interested reader to \cite{krajicek1995bounded} and \cite{cook2010logical} which are monographs treating bounded arithmetic in depth.

In this work, we are concerned with provability in the theory $\PV_1$, introduced in~\cite{krajicek1991bounded} as a first-order extension of Cook's theory $\PV$~\cite{cook1975}, which can be understood in a well-defined sense as a theory of arithmetic with induction accepted only for polynomial-time predicates. Its language contains a function symbol for each polynomial-time algorithm. It is known that the theory $\PV_1$ and its extensions prove many fundamental results of algebra, complexity theory and number theory~\cite{jerabek2005weak,muller2020feasibly,jerabek2022iterated,emildual,EmilAbelian,gaysin2024proof,oliveira2025}. In~\cite{EmilAbelian} it was first observed that the theory~$S^1_2$ proves that every number has a prime divisor: taking the maximal $k$ such that $x$ can be written as a product of $k$ numbers greater than $1$ yields the prime factorization of $x$. The logical principle underlying this argument $\Sigma^b_1\text{-LMAX}$ gives the axiomatization of $S^1_2$ over the base theory $\PV_1$, and is not available in $\PV_1$ itself unless $\text{NP}\subseteq \text{P}/poly$~\cite{krajicek1991bounded}. In this work we show, assuming the hypothesis that non-uniform polynomial-time algorithms cannot factor a constant fraction of all multiples of two $n$-bit primes, for any $n\in \NN$, that $\PV_1$ does not prove that every number has a prime divisor. That is, the $\PV$-sentence \PrimeFactor:
\[(\forall x\geq 2)(\exists y\leq x)(2\leq y \land y\mid z \land (\forall z\leq y)(z\mid y\to (z =1 \lor z= y)))\]
cannot be proved in the theory $\PV_1$ even when extended by the sharply bounded choice scheme $BB(\Sigma^b_0)$, which is also known by the name sharply bounded collection. 

Regarding the plausibility of the assumption, the average-case hardness of factorization against non-uniform polynomial-time adversaries is a well established cryptographic assumption (see~\cite[Section 2.2.4]{goldreich2001foundations}). Moreover, our argument allows for uniform sampling from any finite set of primes (see Theorem~\ref{thrmmain}). In theory, this allows us to instead assume the average-case hardness of factoring where the inputs are taken as products of some subset of primes which are expected to be hard to factorize.

Our main technical tool is the KPT witnessing theorem of~\cite{krajicek1991bounded} which extracts from the provability of $\forall\exists \forall$-sentences in $\PV_1$ a~polynomial-time algorithm which outputs candidates for a witness of the existential quantifier, and each time the algorithm fails it obtains a counterexample to the correctness of the candidate. Such an algorithm is guaranteed to find a correct witness after obtaining a constant number of counterexamples. For the theory $\PV_1+BB(\Sigma^b_0)$, we use an extension of the KPT witnessing theorem due to~\cite{cookthapen2006}, where the algorithm is allowed to output up to polynomially many candidates at once, but is still guaranteed to output a correct one in constantly many steps. To analyze the interactive computations which arise in the witnessing theorems we also use several observations about fields of sets, which are classes of subsets of a given set closed under finite intersections, finite unions and complements.

The work is organized as follows, in Section~\ref{secpreli} we briefly recall basic facts about bounded arithmetic and complexity theory, in Section~\ref{secpv} (Corollary~\ref{crllpv}) we prove the unprovability result for $\PV_1$ and finally in Section~\ref{secpvbb} we extend this result to unprovability in $\PV_1+BB(\Sigma^b_0)$ (Theorem~\ref{thrmpvbb}).

\section{Preliminaries}\label{secpreli}

We assume that the reader is acquainted with basic notions of complexity theory and first-order logic. We use the notation $\floor{-}$ for the floor function and the notation $\abs{-}$ for the binary length of a number. We will refrain from giving a formal definition of the theory $\PV_1$ which can be found in~\cite{krajicek1991bounded,krajicek1995bounded}, instead we will define the true universal theory $T_{\PV}$, which is a proper extension of $\PV_1$. All of our results hold even if we replace $\PV_1$ by $T_{\PV}$, because the presence of true universal sentences does not affect the witnessing theorems. A function computed by a Turing machine is usually understood as a function on binary strings, but in bounded arithmetic we usually understand it as a function on numbers.

\begin{defi}
    Let $\PV$ be a language containing a function symbol $f_M$ for every polynomial-time clocked machine $M$, the intended interpretation of each symbol in $\PV$ is then the function computed by the corresponding machine. The theory $T_{\PV}$ is then axiomatized by the set
    \[\{\varphi\text{ is a universal $\PV$-sentence}; \:\NN\models \varphi\}.\]
\end{defi}

The language of the theory $\PV_1$ indeed satisfies the condition in the previous definition, the readers not familiar with the definition of $\PV_1$ can simply assume the language of $T_{\PV}$ is the minimal language satisfying the condition. We will from now on use $\PV$ to denote the language of $\PV_1$ and we will use the term $\PV$-symbol to either mean a function symbol in $\PV$ or the predicate $g(x_1,\dots,x_k) = 1$, where $g$ is a function symbol in $\PV$.

\subsection{Student-teacher protocols}

In this section, we will recall a formal definition of student-teacher protocols, which are sometimes called counterexample computations~\cite{krajivcek1992no}. This notion is usually left undefined and is simply used as a figure of speech, but we will manipulate these protocols in a non-trivial way and having a formal definition helps the clarity of the arguments.

\begin{defi}
    A \emph{student-teacher protocol} (or just a \emph{protocol}) is a triple $(s,t,c)$, where $c\geq 1$ and $s,t:\NN\to \NN$. 

    Given a protocol $(s,t,c)$ and a number $x$ we define the computation of the protocol $(s,t,c)$ on the input $x$ as the $(2c-1)$-tuple:
    \[(y_1,z_1,y_2, z_2, \dots, y_{c-2},z_{c-2},y_{c-1},z_{c-1},y_c)\]
    where 
    \begin{align*}
        y_1 &= s(x),\\
        z_1 &= t(x,y_1), \\
        y_i &= s(x,z_1,\dots,z_{i-1}), \text{ for }1<i\leq c,\\
        z_i &= t(x,y_1,\dots, y_{i}), \text{ for }1<i<c.
    \end{align*}
    We shall sometimes call the function $s$ \emph{the student}, the function $t$ \emph{the teacher}, the tuple $(y_1,\dots, y_c)$ \emph{the student's answers} and the tuple $(z_1,\dots, z_{c-1})$ \emph{the teacher's replies}.
\end{defi}

\begin{defi}
    Let $\varphi(x,y,z)$ be an open $\PV$ formula. We say a function $t$ is a $\varphi$-correcting teacher on the input $x$ if for every function $s$, every $c\in \NN$ the following is satisfied: The computation of $(s,t,c)$ on input $x$ satisfies for all $i<c$ that \[\text{if }\NN \models \lnot (\forall z)\varphi(x,y_i,z)\text{, then }\NN \models \lnot\varphi(x,y_i,z_i).\]
\end{defi}

The following theorem is the main technical tool underlying our unprovability result, we also include a rephrasing using the language of student-teacher protocols.

\begin{thm}[The KPT theorem~\cite{krajicek1991bounded}]\label{thrmkpt}
    Let $\varphi(x,y,z)$ be an open formula. If
    \[T_{\PV} \vdash (\forall x)(\exists y)(\forall z)(\varphi(x,y,z)),\]
    then there is a number $c\in\NN$, and $\PV$-symbols $f_1,\dots,f_c$ such that
    \[T_{\PV} \vdash \varphi(x,f_1(x),z_1) \lor \varphi(x,f_2(x,z_1),z_2) \lor \dots \lor \varphi(x,f_c(x,z_1,\dots,z_{c-1}),z_c).\]

    Moreover, there is a polynomial-time function $s$, such that for any input $x$ and any teacher $t$ which is $\varphi$-correcting on the input $x$, we have that the computation of the protocol $(s,t,c)$ on an input $x$ contains some $y_i$ which satisfies $\NN \models (\forall z)(\varphi(x,y_i,z))$. Note that in general the running time of $s(x,z_1,\dots,z_i)$ is a multivariate polynomial in $\abs{x}, \abs{z_1},\dots,\abs{z_i}$, but in this work each value $\abs{z_i}$ is always polynomial in $\abs{x}$.
\end{thm}

\section{The unprovability in $\PV_1$}\label{secpv}

The following sentence $\PrimeFactor$ formalizes the statement `every number has a prime factor'. In the rest of this section, we will establish unprovability of this sentence in the theory $T_{\PV}$.

\begin{defi}
    Let $\PrimeFactor_0(x,y,z)$ be the $\PV$-formula
    \[(2\leq  y) \land (y\mid x) \land (z\mid y \to (z=1 \lor z=y)),\]
    where $a \mid b$ denotes the $\PV$-symbol for the divisibility relation.
    
    Moreover, let $\PrimeFactor$ be the $\PV$-sentence
    \[(\forall x\geq 2)(\exists y\leq x)(\forall z\leq y)(\PrimeFactor_0(x,y,z)).\]
\end{defi}

The main idea behind the unprovability is using the KPT theorem, obtaining the polynomial-time student $s$ and bringing its existence to a contradiction with the assumption about the hardness of factoring. Our goal is to prove that there is a teacher which can be simulated in polynomial time and forces the student to do some non-trivial factorization (at least for a large fraction of input parameters).

\begin{defi}\label{defiteacher}
    Assume $s$ is a polynomial-time function, $c,d\geq 1$, and $p_1,\dots,p_d$ are distinct primes. We define a function $t_{(p_1,\dots,p_{d})}$ which serves as a teacher in the protocol $(s,t_{(p_1,\dots,p_d)},c)$ on the input $x=\prod_{i=1}^d p_i$ as follows:
    \begin{enumerate}
        \item (Student's answers) If the student's last answer was not a divisor of $x$ which is greater than $1$, then the $t_{(p_1,\dots,p_d)}$ simply outputs $1$. We will from now on assume the student's answers are always divisors of $x$ which are greater than $1$, as we have already defined the teacher's behavior on the other answers.
        \item (Obvious numbers) Let $1\leq i\leq c$ and assume the student's answers $y_1,\dots,y_{i-1}$ are given, and teacher's replies $z_1,\dots,z_{i-1}$ are also given. We say a number is \emph{obvious (at round $i$)} if it can be obtained from the set $S=\{x,y_1,\dots,y_{i-1},z_1,\dots,z_{i-1}\}$ by $\gcd$ and division without remainder. That is, the prime factorization of an obvious number can be obtained from the prime factorizations of the numbers in $S$ by unions, intersections and complements. A prime factorization of an obvious number is called an \emph{obvious set (at round $i$)}. Note, that the only obvious numbers at round~$1$ are $1$ and $x$.
        \item (Teacher's replies) Let $1\leq i\leq c$. Assume the student's answers $y_1,\dots,y_{i}$ are given, and teacher's replies $z_1,\dots,z_{i-1}$ are also given. The teacher's reply $z_i$ is then one of the following:
        \begin{enumerate}
            \item If $y_i = p_j$ for some $1\leq j \leq d$, then $z_i=p_j$.
            \item Otherwise, if the $\gcd$ of $y_i$ and some obvious number is a proper divisor of $y_i$, then output smallest such $\gcd$.
            \item If neither (a) nor (b) hold, assume that the prime factorization of  $y_i$ is $p_{i_1},\dots,p_{i_{l}}$, $2\leq l\leq d$ and $1\leq i_j\leq d$ for every $j\in \{1,\dots,l\}$. Then, we put $z_i=p_{i_1}\cdot \cdots p_{i_{\floor{l/2}}}$, in which case we say that the teacher divided the student's answer by every value $p_{i_j}$ at round $i$, where $\floor{l/2} < j\leq l$.
        \end{enumerate}
    \end{enumerate}
\end{defi}

Note that the teacher $t_{(p_1,\dots,p_d)}$ is always $\PrimeFactor_0$-correcting on the input $x=\prod_{i=1}^d p_i$. Moreover, if we fix a student $s$ and assume that at the first $i$ rounds, $1\leq i \leq c-1$, the teacher $t_{(p_1,\dots,p_d)}$ divides only by primes from some set $\{p_{i_1},\dots,p_{i_k}\}$, then for $j\in\{1,\dots,i\}$ the reply $z_j$ can be computed by a polynomial-time algorithm which has access to $x,y_1,\dots,y_j,p_{i_1},\dots,p_{i_k}$, where $y_1,\dots,y_j$ are the answers of $s$.

\begin{defi}
    Assume $s$ is a polynomial-time function, $c,d\geq 1$, and $p_1,\dots,p_d$ are distinct primes. Letting $1\leq l <k \leq d$, we say that $s$ with $\{p_1,\dots,p_d\}$ \emph{breaks} $p_lp_k$ if for some permutation $\pi$ on $\{1,\dots,d\}$ there exists $i<c$ such that the computation of protocol $(s,t_{(p_{\pi(1)},\dots, p_{\pi(d)})},c)$ on the input $x=\prod_{i=1}^d p_i$ contains the value $y_i$ which satisfies $\gcd(y_i,p_kp_l)\in\{p_k, p_l\}$ and for every $j<i$ the set of numbers the teacher divided by at round $j$ either contains both $p_l$ and $p_k$ or it contains neither of them.
\end{defi}

To analyze the protocol between a student $s$ and the teacher $t_{(p_1,\dots,p_d)}$ we defined, we will need to analyze the system of obvious sets at a given round, which forms a structure called a field of sets. We will recall its definition, and prove two lemmas about it, which will be used later.

\begin{defi}
    A field of sets $\F$ over $X$ is a family of subsets of $X$ closed under finite union, finite intersection and complements. For $S \subseteq \P(X)$ we define the field of sets generated by $S$, denoted $\C(S)$, as the set of all sets which can be obtained from elements of $S$ by iterated application of finite union, finite intersection and complements. An atom in a field of sets $\F$ is an element which is non-empty and none of its non-empty subsets is in $\F$.
\end{defi}

\begin{lem}\label{lemmdiv}
    Let $\F$ be a field of sets, let $A\in \F$ be an atom of $\F$. Let $A'$ be a proper non-empty subset of $A$, then every atom of $\C(\F\cup \{A'\})$ is already an atom in $\F$ or either one of $A'$ and $A\setminus A'$.
\end{lem}
\begin{proof}
    Let $\F$ be a field of sets over $X$ and let $A\in \F$ be an atom. Let us fix some $B\in \C(\F \cup \{A'\})$, thus there is a finite subset $\{A_1,\dots,A_k\}\subseteq \F$ such that $B$ can be written as a term using the operation of conjunction, disjunction and complement using $\{A_1,\dots,A_k\}$ and $A'$. By the De Morgan laws and the~mutual distributivity of union and intersection, there is $C\in \F$ and indices $l_{i,j}$, $k_{i,j}$ such that
    \[B = \bigcup_i \bigcap_j (A_{l_{i,j}}\cap A') \cup \bigcup_i \bigcap_j (A_{k_{i,j}}\cap (X \setminus A')) \cup C.\]

    Since $A$ is an atom in $\F$ and $A'\subseteq A$, then for every $i$ we have that \[\bigcap_j (A_{l_{i,j}}\cap A')\in \{\varnothing, A'\}.\]
    Moreover, for every $i$ and $j$, we either have that $A \subseteq A_{k_{i,j}}$ and thus \[A_{k_{i,j}}\cap (X\setminus A') = A\setminus A',\]
    or $A\cap A_{k_{i,j}} = \varnothing$ and thus $A_{k_{i,j}}\cap (X\setminus A')=A_{k_{i,j}}$.

    This implies that there are $D\in\{\varnothing, A'\}$, $E\in \{\varnothing,A\setminus A'\}$ and $F\in \F$, such that $B = D \cup E \cup F$. Any combination of choices for $D$, $E$ and $F$ then implies that either $B$ is not an atom, or if it is, then either $B\in\{A',A\setminus A'\}$ or $B$ was already an atom in $\F$.
\end{proof}

\begin{lem}\label{lemmbreaks}
    Let $\F$ be a field of sets over $X$ and let $A\subseteq X$ and $A\not \in \F$. Then there are distinct $a,b\in X$ satisfying $\abs{A\cap \{a,b\}}=1$, such that for every $B \in \F$ we have $\abs{B\cap \{a,b\}} \in \{0,2\}$.
\end{lem}
\begin{proof}
    Assume, that for every distinct $a,b\in X$ which satisfy $\abs{A \cap \{a,b\}} = 1$ there is some $B\in \F$ satisfying $\abs{B\cap \{a,b\}} = 1$. This along with $\F$ being closed under complements implies that for every $a\in A$ and $b\in X\setminus A$ there is a $B\in \F$ such that $\{a,b\} \cap B = \{a\}$.

    Define for every $a\in A$ and $b\in X\setminus A$ the set $B^{a,b}\in \F$ as the set $B$ from the previous sentence. Then, for every such $a$ and $b$ we have
    \[\{a\}\subseteq \bigcap_{b\in X\setminus A} B^{a,b} \subseteq A,\]
    thus $A = \bigcup_{a\in A} \bigcap_{b\in X\setminus A} B^{a,b},$ a contradiction.
\end{proof}

\begin{lem}\label{lemmaobvsize}
    Assume $s$ is a polynomial-time function, $c\geq 1$, $d=2^c$, $p_1,\dots,p_d$ are distinct primes and $x=\prod_{i=1}^d p_i$.  Assume that there is an $i$ such that for all $j\leq i$ we have that $y_j=s(x,z_1,\dots,z_{j-1})$ is a number which is obvious at round $j$ during the interaction with the teacher $t_{(p_1,\dots,p_d)}$. Then the number of distinct prime factors of an obvious number at round $i$ is either $0$ or at least $2^{c-i+1}$.
\end{lem}
\begin{proof}
    By induction on $i$ we will prove that the size of every minimal non-empty obvious set is a power of two which is at least $2^{c-i+1}$. For $i=1$ the only obvious number with a non-empty prime factorization is $x$ itself with $2^{c}$-many prime factors. Assume the statement holds for $i$, the answer of $s$ is $y_i$, which is an obvious number at round $i$, and the teacher $t_{(p_1,\dots,p_d)}$ replies with $z_i$. If $z_i$ was obvious at round $i$, no new obvious numbers are introduced at round $i+1$. 
    
    Assume that $z_i$ was not obvious at round $i$. It is straightforward to check that the set of all obvious sets at round $i$ forms a field of sets, with atoms corresponding to prime factorizations of obvious numbers which are minimal with respect to division. By Lemma~\ref{lemmdiv}, where $\F$ is taken to be the field of sets at round $i$ and $A'$ is taken as the prime factorization of $z_i$, we obtain that the only new atoms in the field of obvious sets at round $i+1$ are prime factorizations of $z_i$ and $y_i/z_i$, both of size $2^{c-(i+1)+1}$. Since the cardinalities of sets in the field of obvious sets at round $i$ correspond to the number of distinct prime factors of obvious numbers, we obtain that every obvious number at round $i+1$ has either $0$ prime divisors, or the prime factorization contains an atom as a subset, which contains at least least $2^{c-(i+1)+1}$ distinct prime divisors.
\end{proof}

\begin{lem}\label{lemmabreak}
    Assume $s$ is a polynomial-time function, $c\geq 1$, $d=2^c$ and $p_1,\dots,p_d$ are distinct primes such that for any teacher $t$ which is $\PrimeFactor_0$-correcting on the input $x=\prod_{i=1}^d x_i$, the computation of $(s,t,c)$ on the input $x$ contains some prime factor of $x$ as one of the student's answers $y_i$.
    
    Then, there are distinct indices $l,k\in \{1,\dots,d\}$, such that $s$ with $\{p_1,\dots,p_d\}$ breaks $p_lp_k$. 
\end{lem}
\begin{proof}
    We will analyze the computation of the protocol $(s,t_{(p_1,\dots,p_{d})},c)$ on the input $x=\prod_{i=1}^{d}p_i$. Assume such a value exists, consider the smallest $1\leq i\leq c$ such that the value $y_i = s(x,z_1,\dots,z_{i-1})$ is not obvious at round $i$, and let $A$ be the set of prime factors of $y_i$. The set of all obvious sets at round $i$ forms a field of sets which does not contain $A$, therefore by Lemma~\ref{lemmbreaks} there are distinct $l,k \in \{1,\dots,d\}$ such that $\gcd(p_lp_k, y_i)\in\{p_l,p_k\}$ and the teacher did not divide any answer of the student by exactly one of $p_l$ and $p_k$, thus $s$ with $\{p_1,\dots,p_d\}$ breaks $p_lp_k$. 
    
    By Lemma~\ref{lemmaobvsize}, if the Student $s$ does not respond with a non-obvious number then the number of prime factors of the value $y_c = s(x,z_1,\dots,z_{c-1})$ is at least $2$, a contradiction with the assumption on $s$.
\end{proof}

\begin{lem}\label{lemmprob}
    Let $\Omega$ be a set, $d\geq 1$ and $F$ a function from subsets of $\Omega$ of size $d$ such that
    \[\forall T\subseteq \Omega, \abs{T}=d: F(T)\text{ is a non-empty set of $2$-element subsets of $T$.} \]
    Then,
    \[\Pr_{x_1,\dots,x_d \sim \Omega}[\{x_1,x_2\} \in F(\{x_1,\dots,x_d\})|x_1,\dots,x_d \text{ are distinct}]\geq \binom{d}{2}^{-1},\]
    where $x_1,\dots,x_d$ are sampled uniformly and independently.
\end{lem}
\begin{proof}
    Let $T$ be a fixed $d$-element subset of $\Omega$, and let $\{y,z\}\in F(T)$. Random choice of $x_1,\dots,x_d$ such that $T=\{x_1,\dots,x_d\}$ will satisfy $\{x_1,x_2\}=\{y,z\}$ with probability $\binom{d}{2}^{-1}$. For different choices of $T$ the events \[E_T=\{(x_1,\dots,x_d) ; T=\{x_1,\dots,x_d\}\}\] are disjoint, and thus the statement of the lemma follows.
\end{proof}

\begin{thm}\label{thrmmain}
    Assume $s$ is a polynomial-time function and $c,d\geq 1$ such that for any distinct primes $p_1,\dots,p_d$, there are some $1\leq l <k \leq d$ such that $s$ with $\{p_1,\dots,p_d\}$ breaks $p_lp_k$ (during a $c$-round computation). 
    Then there is $r>0$ and a polynomial time function $f$ such that for every finite set of primes $D$ we have
    \[\Pr_{\substack{p,q \sim D\\ p_1,\dots,p_{d-2} \sim D}}  \left[f(pq,p_1,\dots,p_{d-2}) \in \{p,q\}\right]\geq r,\]
    where $p,q,p_1,\dots,p_{d-2}$ are sampled uniformly and independently.
\end{thm}
\begin{proof}

    Consider the following algorithm $f$: On the input $(pq,p_1,\dots,p_{d-2})$, it first checks whether all numbers $p_1,\dots,p_{d-2}$ are distinct and do not divide $pq$ and then it tries to simulate for every permutation $\pi$ on $\{p,q,p_1,\dots,p_{d-2}\}$ protocols between the student $s$ and the teacher $t_{(\pi(p),\pi(q),\pi(p_1),\dots,\pi(p_{d-2}))}$ on the input $x=pq\prod_{i=1}^{d-2}p_i$ in the following way:
        
        It iterates over every permutation of $\{*_1,*_2, p_1,\dots,p_{d-2}\}$, where $*_1$~and~$*_2$ are placeholder values for $p$ and $q$ which are unknown to $f$. It then simulates the communication with the teacher $t_{(\pi(p),\pi(q),\pi(p_1),\dots,\pi(p_{d-2}))}$, where $\pi$ is the induced permutation, until division by exactly one of $p$ or $q$ is needed to proceed, in which case the simulation is aborted. Importantly, if it happens that the given permutation leads to $s$ with $\{p,q,p_1,\dots,p_{d-2}\}$ breaking $p,q$, then the knowledge of either $p$ or $q$ was not needed in the simulation of the teacher, as the knowledge of the product $pq$ suffices for every answer of the teacher. 
        
        After each answer of the student $s$, the algorithm $f$ tries to take $\gcd$ of it and $pq$ and if it finds a proper divisor it returns it.

    Let $F(p_1',\dots,p_d')$ be a function which outputs the set of all pairs of primes which are broken by the student $s$ with $\{p_1',\dots,p_d'\}$. By Lemma~\ref{lemmabreak}, this set is always non-empty and thus $F$ satisfies the conditions of Lemma~\ref{lemmprob}, which implies 
    \[\Pr_{p,q,p_1,\dots,p_{d-2}\sim D}[f(pq,p_1,\dots,p_{d-2})\in\{p,q\}|p,q,p_1,\dots,p_{d-2}\text{ distinct}]\geq \binom{d}{2}^{-1}.\]

    If $\abs{D}\leq 4 \binom{d}{2}$, then the probability of $p_1$ dividing $pq$ is at least $1/(4\binom{d}{2})$, and if $\abs{D}\geq 4\binom{d}{2}$, then the probability of $p,q,p_1,\dots,p_{d-2}$ being all distinct is at least \[\prod_{i=0}^{d-1} \left(1-\frac{i}{\abs{D}}\right ) \geq \left(1- \frac{d}{4 \binom{d}{2}}\right)^{d-1} \geq (1/2),\]
    and thus, combining this with the previous paragraph the probability that $f$ finds a factor of $pq$ given $p_1,\dots, p_{d-2}$ is at least $\frac{1}{4\binom{d}{2}}$.
\end{proof}

\begin{cor}\label{crllpv}
    Assume that for every $r>0$ and every sequence of Boolean circuits $\{C_n\}_{n\in\NN}$ of polynomial size there is an $n$ such that
    \[\Pr_{p,q\sim P_n}[C_n(pq) \in \{p,q\}] < r,\]
    where $P_n$ is the set of all primes of length $n$ and $p,q$ are sampled uniformly and independently. Then,
    $T_{\PV} \nvdash \PrimeFactor.$
\end{cor}
\begin{proof}
    We will prove the contrapositive, therefore we assume that \[T_{\PV}\vdash \PrimeFactor.\] By Theorem~\ref{thrmkpt} there is a polynomial-time function $s$ and $c\in \NN$ satisfying the assumptions of Lemma~\ref{lemmabreak} and therefore also Theorem~\ref{thrmmain}, which implies that there is $r>0$, $d\in \NN$ and a polynomial time function $f$, which for every $D=P_n$ satisfies
        \[\Pr_{\substack{p,q \sim D\\ p_1,\dots,p_{d-2} \sim D}}  \left[f(pq,p_1,\dots,p_{d-2}) \in \{p,q\}\right]\geq r.\]
    By an averaging argument, this means that for each $n$ there are specific elements $p_1,\dots,p_{d-2}\in P_n$ such that 
    \[\Pr_{p,q\sim P_n}[f(pq,p_1,\dots,p_{d-2}) \in \{p,q\}]\geq r,\]
    and thus we can take $C_n$ to be the circuit computing $f(-,p_1,\dots,p_{d-2})$.
\end{proof}

\section{The Unprovability with the Sharply Bounded Choice Scheme}\label{secpvbb}

We will establish the unprovability in the theory $T_{\PV}+BB(\Sigma^b_0)$, which extends $T_{\PV}$ by the following scheme.

\begin{defi}
    The sharply bounded choice scheme $BB(\Sigma^b_0)$ is the set of axioms of the form
    \[(\forall i\leq \abs{a})(\exists y\leq a)(\varphi(i,y)) \to (\exists w)(\forall i \leq \abs{a})(\varphi(i,[w]_i)),\]
    for each $\varphi\in \Sigma^b_0$, where $[w]_i$ is a shorthand for a $\PV$-symbol which outputs the $i$-th member of the sequence coded by $w$.
\end{defi}

We also need the following variant of the KPT theorem for $T_{\PV}+BB(\Sigma^b_0)$.

\begin{thmC}[\cite{cookthapen2006}]\label{thrmkptbb}
    Let $\varphi(x,y,z)$ be an open formula. If
    \[T_{\PV} + BB(\Sigma^b_0) \vdash (\forall x)(\exists y)(\forall z)(\varphi(x,y,z)),\]
    then there is a number $c\in\NN$, and $\PV$-symbols $b,f_1,\dots,f_c$ such that
    \begin{align*}
        T_{\PV} \vdash (\exists i\leq \abs{b(x)})&\varphi(x,[f_1(x)]_i,[z_1]_i)\\ \lor (\exists i \leq \abs{b(x)})&\varphi(x,[f_2(x,z_1)]_i,[z_2]_i) \\ & \vdots  \\ \lor (\exists i\leq \abs{b(x)})&\varphi(x,[f_c(x,z_1,\dots,z_{c-1})]_i,[z_j]_i).
    \end{align*}
\end{thmC}

Note that even though the formula $(\exists i \leq \abs{b(x)})\varphi(x,[y]_i,[z]_i)$ is not open, it is actually equivalent to an open formula as the existential quantifier is sharply bounded, that is, the bound's outermost function symbol is $\abs{-}$. A $\PV$-symbol $g$ which tries all possible values for $i$ and outputs $1$ if and only if the open kernel is satisfied for at least one $i$ is straightforward to construct, and so the formula is then equivalent in $T_{\PV}$ (even in $\PV_1$) to $g(x,y,z)=1$, which is an open formula. We will call this formula $\PrimeFactor_1$.

We will now define a parallel variant of the teacher $t_{(p_1,\dots,p_d)}$, which replies to a student outputting sequences of divisors of $x$.

\begin{defi}\label{defipteacher}
    Assume $s^P$, $b$ are polynomial-time functions, ${e}\geq 1$, $d=2^{{e}}$ and $p_1,\dots,p_d$ are distinct primes. We define a function $t^P_{(p_1,\dots,p_{d})}$ which serves as a teacher in the protocol $(s^P,t^P_{(p_1,\dots,p_d)},{e})$ on the input $x=\prod_{i=1}^{d} p_i$ as follows:
    \begin{enumerate}
        \item (Student's answers) The student's answers are interpreted as sequences of divisors of $x$ which are greater than $1$ of length $b(x)$. We will define the answer of the teacher coordinate-wise. That is for every divisor $y_i^j$, where $1\leq i\leq d$ and $j\leq \abs{b(x)}$, we will define the teacher's reply $z^j_i$. In the case the student's answer is not a sequence, the teacher replies with $1$ and if it is a sequence but any element of the sequence is not a divisor of $x$ greater than $1$, the teacher's reply on that coordinate is $1$. 
         \item (Obvious numbers) Let $1\leq i\leq {e}$ and assume the student's answers $y_1,\dots,y_{i-1}$ are given, and teacher's replies $z_1,\dots,z_{i-1}$ are also given. We say a number is \emph{obvious (at round $i$)} if it can be obtained from the set 
         \[S=\{x\}\cup \{y_k^j; j \leq \abs{b(x)}, 1 \leq k \leq i-1\} \cup \{z_k^j; j \leq \abs{b(x)}, 1\leq k\leq i-1\}\] by $\gcd$ and division without remainder. That is, the prime factorization of an obvious number can be obtained from the prime factorizations of the numbers in $S$ by unions, intersections and complements. A prime factorization of an obvious number is called an \emph{obvious set (at round $i$)}. Note, that the only obvious number at round~$1$ are $1$ and $x$.
        \item (Teacher's replies) Let $1\leq i\leq {e}$. Assume the student's answers $y_1,\dots,y_{i}$ are given, and teacher's replies $z_1,\dots,z_{i-1}$ are also given. Let $j\leq \abs{b(x)}$, the teacher's reply on the $j$-th coordinate $z_i^j$ is then one of the following:
        \begin{enumerate}
            \item If $y_i^j = p_k$ for some $1\leq k \leq d$, then $z_i^j=p_k$.
            \item Otherwise, if the $\gcd$ of $y_i^j$ and some obvious number is a proper divisor of $y_i^j$, then output the least such $\gcd$.
            \item Otherwise, assume that the prime factorization of $y_i^j$ is $p_{i_1},\dots,p_{i_{l}}$, $2\leq l\leq d$ and $1\leq i_k\leq d$ for every $k\in \{1,\dots,l\}$. Then, we put $z_i^j=p_{i_1}\cdot \cdots p_{i_{\floor{l/2}}}$, in which case we say that the teacher divided the student's answer by every value $p_{i_j}$ at round $i$, where $\floor{l/2} < j\leq l$.
        \end{enumerate}
    \end{enumerate}
\end{defi}

Note that the teacher $t^P_{(p_1,\dots,p_d)}$ is always $\PrimeFactor_1$-correcting on the input $x=\prod_{i=1}^d p_i$. Moreover, when given access to the primes it divides by at the first $i$ rounds, the replies $z_1,\dots,z_i$ of $t_{(p_1,\dots,p_d)}^P$ on the input $x$ can be computed in polynomial-time. 

We now show that the size of atoms in the field of obvious sets does not shrink too quickly and that the cardinality of all atoms also does not increase too quickly, unless the student outputs some non-obvious number.

\begin{lem}\label{lemmatoms}
    Assume $s^P$, $b$ are polynomial-time functions, ${e}\geq 1$, $d=2^{{e}}$ and $p_1,\dots,p_d$ are distinct primes. Assume that there is an $i$ such that for all $k\leq i$ we have that for\linebreak every $j\leq \abs{b(x)}$: 
    \[y_k^j=[s(P,z_1,\dots,z_{j-1})]_j\text{ is a number which is obvious at round $j$.}\] 
    Then the number of distinct atoms in the field of obvious sets at round $i$ is at most $2^i$, and each of those atoms is of size $2^{c-i+1}$.
\end{lem}
\begin{proof}
    By induction on $i$. For $i=1$, the only atom is the prime factorization of $x=\prod^d_{k=1} p_k$ of size $2^{c}$.

    Assume the statement holds for $i$. Then the student's answer contains at most $2^i$-many atoms at round $i$. The teacher $t^P_{(p_1,\dots,p_d)}$ then replies with a sequence of numbers whose prime factorizations are subsets of the prime factorizations of the student's numbers. By iterated application of Lemma~\ref{lemmdiv}, and the fact 
    \[\C(\F \cup \{A_1,\dots,A_m\})=\C(\C(\cdots \C(\C(\F \cup \{A_1\})\cup \{A_2\}) \cdots )\cup \{A_m\})\]
    for any field of sets $\F$ over $X$ and subsets $A_1,\dots,A_m \subseteq X$, $m\in \NN$, we have that the field of obvious sets at round $i+1$ contains at most two atoms for each of the atoms in the field of obvious sets at round $i$, and those atoms at round $i+1$ are of size $2^{c-i+1}/2$, which concludes the inductive step.
\end{proof}

To finish the proof we will apply Theorem~\ref{thrmmain}. To do so, we will convert the parallel student $s^P$ into a sequential one $s$ by making $s$ output the obvious answers of $s^P$ sequentially, and if $s^P$ were to answer with a non-obvious number, then this number is taken as the answer of $s$ for the rest of the computation. This increases the number of rounds in the protocol from ${e}$ to $2^{{e}}-1$, which is sufficient for our application as the new number of rounds is still a constant.

\begin{lem}\label{lemmconversion}
    Assume that $T_{\PV}+ BB(\Sigma^b_0) \vdash \PrimeFactor$, then there is a polynomial-time student $s$ and $c,d\geq 1$ such that for any distinct primes $p_1,\dots,p_d$, there exists $1\leq l < k \leq d$ such that $s$ with $\{p_1,\dots,p_d\}$ breaks $p_lp_k$ (during a $c$-round computation).
\end{lem}
\begin{proof}
    By Theorem~\ref{thrmkptbb} we obtain a polynomial-time function $s^P$ and ${{e}\in \NN}$ such that for distinct primes $p_1,\dots, p_d, d=2^{{e}},$ the computation of the protocol $(s^P,t^P_{(p_1,\dots,p_d)},{e})$ on the input $x=\prod_{i=1}^d p_i$ contains the student's answer~$y_i$ which is a sequence of length $\abs{b(x)}$ containing at least one prime divisor of~$x$.

    Consider a polynomial-time function $s$ which serves as a student in the protocol $P_0 =(s,t_{(p_1,\dots,p_d)},c)$, where $c=2^{{e}}-1$, which we define as follows: 
    \begin{enumerate}
        \item First, we partition the set $\{1,\dots,2^{{e}}\}$ into the sets $R_i = \{2^{i-1},\dots , 2^{i}-1\},$ where $1\leq i \leq {e}$.
        \item The student $s$ keeps a partial computation of the protocol \[P=(s^P,t^P_{(p_1,\dots,p_d)},{e}).\] At round $2^{i-1}$ the student $s$ will have constructed the following part of the computation:
        \[(x,y_1,z_1,\dots,y_{i-1},z_{i-1},y_i).\]
        \item The answers of the student $s$ at rounds contained in $R_i$, $1\leq i \leq {e}$, are one of the following:
        \begin{enumerate}
            \item If $i=1$, then the student runs $s^P(x)$ and obtains $y_1$, if the output contains at least one non-obvious number, then $s$ outputs it for all of the remaining rounds. Otherwise, it outputs $x$.
            \item If $1<i\leq {e}$, and all of the numbers contained in the answers of $s^P$ have been obvious, then by Lemma~\ref{lemmatoms}, there are at most $2^{i-1}$ atoms in the field of obvious sets at round $i-1$ in $P$, and thus the replies of the teacher $t_{(p_1,\dots,p_d)}$ at rounds in $R_{i-1}$ can be collected to obtain the reply of $t^P_{(p_1,\dots,p_d)}$ which we denote $z_{i-1}$, this in turn allows us to compute the $i$-th reply of $s^P$ which we denote $y_i$. 
            \begin{itemize}
                \item If all numbers in $y_i$ are obvious, then by Lemma~\ref{lemmabreak} there is at most $2^i$ of them, and we use them as the answers of $s$ at rounds in $R_i$ (in any particular order).
                \item If there is a number in $y_i$ which is non-obvious, the student $s$ outputs it for all of the remaining rounds.
            \end{itemize} 
        \end{enumerate}
    \end{enumerate}
    Note that the term `obvious number' is used here in the sense of Definition~\ref{defipteacher} for the computation of the protocol $P$.

    By Lemma~\ref{lemmatoms}, we know that on the input $x$ the student $s^P$ has to output a non-obvious number at some round, otherwise all elements in $y_{{e}}$ have at least two prime divisors contradicting Theorem~\ref{thrmkptbb}. Moreover, if a non-obvious number $a$ is contained in some answer $y_i$ of $s^P$, but all previous answers contained only obvious answers, then the answers of $s$ at rounds in $R_i$ are all $a$, which is non-obvious in the sense of Definition~\ref{defiteacher}. As in the Lemma~\ref{lemmabreak}, this implies that there are distinct $k,l\in \{1,\dots,d\}$ such that $s$ with $\{p_1,\dots,p_d\}$ breaks $p_lp_k$.
\end{proof}

\begin{thm}\label{thrmpvbb}
    Assume that for every $r>0$ and every sequence of Boolean circuits $\{C_n\}_{n\in\NN}$ of polynomial size there is an $n$ such that
    \[\Pr_{p,q\sim P_n}[C_n(pq) \in \{p,q\}] < r,\]
    where $P_n$ is the set of all primes of length $n$ and it is sampled uniformly and independently. Then,
    $T_{\PV} + BB(\Sigma^b_0) \nvdash \PrimeFactor.$
\end{thm}
\begin{proof}
    The proof mirrors the one of Corollary~\ref{crllpv}. We assume that the theory $T_{PV}+BB(\Sigma^b_0)$ actually does prove $\PrimeFactor.$ By Lemma~\ref{lemmconversion} there is a polynomial-time function $s$ and $c\geq 1$ satisfying the assumptions of Theorem~\ref{thrmmain}, which if we combine with an averaging argument gives us the sequence of circuits with desired properties.
\end{proof}

Let us note that the idea which we applied in this section of unfolding the query-sequences into constantly many individual queries uses the fact that the number of obvious sets which are atoms at round $i$ is a function of $i$ and does not depend on the length of $x$. For example, the argument fails if we were to try to reduce to provability of the sharply bounded choice scheme itself into a~constant number of single queries, as the number of possible indicies which the student can query in the round $i$ is $\abs{b(x)}$, which does depend on the length of $x$.

\section{Concluding remarks}\label{secconcl}

The unprovability of the formula $\PrimeFactor$ in $T_{\PV}+ BB(\Sigma^b_0)$ implies the existence of a model
$M \models T_{\PV} + BB(\Sigma^b_0)+\lnot \PrimeFactor$. In other words, there is $m\in M$ such that
\[M\models (\forall y)(\exists z)((y \neq 1 \land y\mid m) \to (z\mid y \land z\neq 1 \land z\neq y)),\]
meaning that every divisor of $m$ has a proper divisor in $M$. This can be rephrased as saying that $m$ has no irreducible divisors, which also implies that $m$ has no prime factorization in $M$. \emph{A Furstenberg domain}~\cite{clark2017} is an integral domain in which every non-invertible element has an irreducible divisor. The main result of this paper can thus be restated as: There are models of $T_{\PV}+BB(\Sigma^b_0)$ which are not positive parts of a Furstenberg domain. Let us recall that every model of $S^1_2$ is indeed a positive part of a Furstenberg domain.

In~\cite{JERABEK2016380} Jeřábek has shown that integer factorization search problem is probabilistically polynomial-time many-one reducible to the $\textsc{WeakPigeon}$ problem, whose instances are functions $f:\{0,1\}^n\to \{0,1\}^{n-1}$ represented by a Boolean circuit and solutions are pairs $x_1,x_2$ of distinct inputs such that $f(x_1)=f(x_2)$. The reduction admits non-uniform derandomization, and thus our assumption also implies that $\text{PV}_1$ also does not prove the injective weak pigeonhole principle $\text{iWPHP}$, but it is not clear that the two unprovabilities imply the existence of a model where both statements fail. It is also not clear whether one could find an assumption which implies that $\PV_1 +\PrimeFactor\nvdash \text{iWPHP}$, as there could plausibly be a some reduction from $\textsc{WeakPigeon}$ to integer factorization.
\begin{qu}
    Is there a model $M\models \text{PV}_1+\PrimeFactor+\lnot\text{iWPHP}$ under some plausible assumption?
\end{qu}

In~\cite{cookthapen2006} Cook and Thapen have shown, assuming factorization is not in probabilistic polynomial time, that $\PV_1 \nvdash BB(\Sigma^b_0)$. Our assumption on the hardness of factorization is stronger then theirs, and thus after combining our result with theirs we obtain (under this hypothesis):
\[ \PV_1 \lneq \PV_1+ BB(\Sigma^b_0) \lneq S^1_2.\]
This gives a single assumption which can separate three consecutive theories contained in Buss's hierarchy. It would be interesting to see if we can get separations from even higher theories from our assumption.

\begin{qu}
    Assume that there is no polynomial-size family of Boolean circuits which can factorize a constant fraction of all products of two $n$-bit primes for every $n\in\NN$. Can we show that  $S^1_2\lneq T^1_2$?
\end{qu}

Note that assuming $\text{L}^{\text{NP}} \neq \text{P}^{\text{NP}}$ we can separate $S^1_2$ from $T^1_2$~\cite{krajivcek1993fragments} and assuming that the polynomial hierarchy does not collapse, we can also separate $S^i_2$ from $S^{i+1}_2$ for any $i\geq 1$~\cite{krajicek1991bounded}. 

More generally, we can ask about other hypotheses, which do not explicitly mention the complexity classes used in the definition of the theories, but which separate as many consecutive theories of bounded arithmetic as possible. One possible interpretation of the question is the following.

\begin{qu}
    Is there an assumption not explicitly mentioning the polynomial hierarchy which separates $S^i_2$ and $S^{i+1}_2$ for any $i\geq 1$?
\end{qu}

\section*{Acknowledgment}

The author thanks Jan Krajíček for his guidance and helpful comments. Part of this work was completed when the author was hosted by Igor Oliveira at the University of Warwick. The author also thanks Raheleh Jalali, Erfan Khaniki, Mykyta Narusevych, Daria Pavlova and Neil Thapen for helpful discussions.

\bibliographystyle{alphaurl}
\bibliography{main}

\end{document}